\documentclass[12pt,a4paper]{article}
\usepackage{amsmath}

\usepackage{amssymb,latexsym}

\newtheorem{lemma}{Lemma} [section]

\title{Extremal eigenvalues of the Laplacian in a conformal class of metrics~:
the "conformal spectrum" }
\author{Bruno Colbois and Ahmad El Soufi}
\date{}

\begin{document}

\maketitle

\bigskip
\maketitle

\bigskip

\begin{abstract}

Let $M$ be a compact connected manifold of dimension $n$ endowed with
a conformal class $C$ of
Riemannian metrics of volume one. For any integer $k\geq0$, we
consider the conformal invariant $\lambda_k ^c (C)$ defined as
the supremum of the $k$-th eigenvalue $\lambda_k (g)$ of the
Laplace-Beltrami operator $\Delta_g$, where $g$ runs over $C$.

First, we give a sharp universal lower bound for $\lambda_k ^c (C)$
extending to all $k$ a result obtained by Friedlander and
Nadirashvili for $k=1$. Then, we show that the sequence $
\left\{\lambda_k ^c (C) \right\}$, that we call "conformal spectrum",
is strictly increasing and satisfies, $\forall k\geq 0$,
$\lambda_{k+1} ^c (C)^{n/2} - \lambda_k ^c (C)^{n/2} \geq n^{n/2}
\omega_n $, where $\omega_n $ is the volume of the $n$-dimensional
standard sphere.

When $M$ is an orientable surface of genus $\gamma$, we also consider
the supremum $\lambda_k ^{top} (\gamma)$ of $\lambda_k(g)$ over the
set of all the area one Riemannian metrics on $M$, and study the
behavior of $\lambda_k ^{top} (\gamma)$ in terms of $\gamma$.

\end{abstract}

\bigskip

{\bf 2000 Mathematics Subject Classification}~: $58 J 50$, $58 E 11$,
$35P15$.

\medskip

{\bf Keywords}: Laplacian, eigenvalue, conformal metric, universal lower bound.

\bigskip
\bigskip
\bigskip
\noindent
B. Colbois~: Universit\'e de Neuch\^atel, Laboratoire de
Math\'ematiques, 13 rue E. Argand,
2007 Neuch\^atel, Switzerland.

\smallskip
\noindent
E-mail: Bruno.Colbois@unine.ch

\bigskip
\noindent
A. El Soufi~: Universit\'e de Tours, Laboratoire de Math\'ematiques
et Physique Th\'eorique, Parc de Grandmont, 37200 Tours, France.

\smallskip
\noindent
E-mail: elsoufi@univ-tours.fr.
\eject

\section{Introduction and statement of results}\label{1}

Let $M$ be a closed connected differentiable manifold of dimension $n
\geq 2$. Given a Riemannian metric $g$ on $M$, let
$$spec(g)= \left\{ 0= \lambda_0(g) < \lambda_1(g) \leq \lambda_2(g)
\cdot \cdot \cdot \leq \lambda_k (g) \leq \cdot \cdot \cdot
\right\}$$ be the spectrum of the Laplace-Beltrami operator defined
by $g$.

One of the main topics in spectral geometry is the study of the
variational properties of the functional
$g \longmapsto \lambda_k(g) $ and the finding out of extremal
geometries for $\lambda_k$. Problems of this kind were first studied
in the setting of Euclidean domains where many Faber-Krahn type
inequalities have been established (see [He] for a recent survey). In
the case of closed manifolds which interests us here, the first
result was obtained by
Hersch [H]: on the 2-sphere $ \mathbb S ^2$, the standard metric
$g_s$ maximizes $\lambda_1$ among all the Riemannian metrics of the
same area. Moreover, $g_s$ is, up to isometry, the unique maximizer.

Recall that the behavior of $\lambda_k$ under scaling of the metric
is given by
$\lambda_k (cg)={\lambda_k (g) / c}$. Hence, a normalization is
required. The metric invariant usually considered for this
normalization is the volume $V(g)$. Therefore,
we denote by ${\cal M}(M)$ the set of all Riemannian metrics of {\it
volume one} on $M$
and, for any $g \in {\cal M}(M)$, we set
$$ [g]= \left\{ g'\in {\cal M}(M) \mid g'{\ } {\hbox {is conformal
to}} {\ }g \right\}.$$

It is well known that if $M$ is of dimension $n \geq 3$, then,
$\forall k \geq1$, $\lambda_k$
is not bounded on ${\cal M} (M)$ (see [CD]). On the other hand,
Korevaar [K] showed that in dimension 2, $\lambda_k$ is bounded on
${\cal M} (M)$ and that, in all dimensions, the restriction of
$\lambda_k$ to
any conformal class of metrics of fixed volume is bounded. Hence, for
any natural integer $k$ and any conformal class of metrics $[g]$ on
$M$, we define the {\it conformal $k$-th eigenvalue} of $(M, [g])$ to
be
$$ \lambda_k ^c (M, [g]) = \sup_{g'\in [g]} \lambda_k (g') = \sup
\left\{ \lambda_k (g') V(g')^{2/n} \mid g'{\ } {\hbox {is conformal
to}} {\ }g \right\}.$$
The sequence $\left\{\lambda_k ^c (M,[g]) \right\}$ constitutes the
{\it conformal spectrum} of $(M,[g])$.

In dimension 2, one can also define a {\it topological spectrum} by
setting, for any genus $\gamma$ and any integer
$ k\geq0$,
$$\lambda_k ^{top} (\gamma) = \sup \left\{ \lambda_k (g) \mid g \in
{\cal M}(M_{\gamma}) \right\},$$
$ M_{\gamma}$ being an orientable compact surface of genus $\gamma$.

 The aim of this paper is to emphasize some properties of the
conformal and topological spectra. Let us first recall some of the
known results. Actually, most of them concern only the first positive
eigenvalue. Indeed, the result of Hersch mentioned above reads:
$\lambda_1 ^{top} (0)=\lambda_1 (g_s)=8\pi$, where $g_s$ is the
standard metric normalized to volume one. For genus one surfaces,
Nadirashvili [N1] showed that $\lambda_1 ^{top} (1)= \lambda_1(g_e)= 8
\pi^2/\sqrt{3}$, where $g_e$ is the flat metric induced on the
2-torus from an equilateral lattice of ${\mathbb R}^2$. For
arbitrary genus, Yang and Yau [YY] proved the following inequality
(see also [EI1]):
$$\lambda_1 ^{top} (\gamma) \leq 8 \pi [{{\gamma + 3} \over 2} ],$$
where [ ] denotes the integer part, and Korevaar [K] showed the
existence of a universal constant $C$ such that, $\forall k\geq0$,
$$\lambda_k ^{top} (\gamma) \leq C(\gamma + 1)k.$$
In higher dimension, Korevaar also obtained in [K] the estimate~:
$$ \lambda_k ^c (M, [g]) \leq C([g]) k^{2/n}$$
for some constant $C([g])$ depending on $n$ and on a lower bound of $Ric
{\ } d^2$, where $Ric$ is the Ricci curvature and $d$ is the diameter
of $g$ or of another representative of $[g]$.

Regarding the conformal first eigenvalue, the second author and Ilias
[EI2] gave a sufficient condition for a Riemannian metric $g$ to
maximize $\lambda_1$ in its conformal class $[g]$: if there exists a
family $f_1, f_2, \cdot \cdot \cdot, f_p$ of first eigenfunctions
satisfying $\sum_{i} df_i \otimes df_i=g$, then $ \lambda_1 ^c (M,
[g])= \lambda_1 (g)$. This condition is fulfilled in particular by
the metric of any homogeneous Riemannian space with irreducible
isotropy representation. For instance, the first conformal
eigenvalues of the rank one symmetric spaces endowed with their
standard conformal classes $[g_s]$, are given by

\begin{itemize}

 \item $\lambda_1 ^c
({\mathbb S}^n, [g_s])=n\omega_n^{2/n}$, where $\omega_n $ is the
volume of the $n$-dimensional Euclidean sphere of radius one,

\item $\lambda_1 ^c ( {\mathbb R}P^n, [g_s])=2^{n-2 \over n}
(n+1)\omega_n^{2/n}$,

\item $\lambda_1 ^c ( {\mathbb C}P^d, [g_s])=4\pi
(d+1) d!^{-1/d}$,

\item $\lambda_1 ^c ( {\mathbb H}P^d, [g_s])=8\pi (d+1)
(2d+1)!^{-1/2d}$,

\item $\lambda_1 ^c ( {\mathbb C}aP^2, [g_s])=48\pi
({6\over 11!})^{1/8}=8\pi \sqrt 6({9\over 385})^{1/8}.$
\end{itemize}

On the other hand, Ilias, Ros and the second author [EIR] proved that
if
$\Gamma={\mathbb Z}e_1+ {\mathbb Z}e_2 \subset {\mathbb R}^2$ is a
lattice such that $|e_1 | =|e_2|$, then the corresponding flat metric
$g_{_\Gamma}$ on ${\mathbb T}^2$ satisfies $ \lambda_1 ^c ({\mathbb
T}^2, [g_{_\Gamma}])=\lambda_1 (g_{_\Gamma})$. A higher dimensional
version of this result was also established in [EI3]. Neverthless,
the authors [CE] have showed that when the length ratio $|e_2| /
|e_1 |$ of the vectors $e_1$ and $e_2$ is sufficiently far from 1,
then $ \lambda_1 ^c ({\mathbb T}^2, [g_{_\Gamma}])> \lambda_1
(g_{_\Gamma})$, that is, $g_{_\Gamma}$ does not maximize $\lambda_1$
on $[g_{_\Gamma}]$.

Finally, the following relationship between $\lambda_1 ^c (M, [g])$
and the conformal volume $V_c(M,[g])$ is due to Li and Yau [LY] in
dimension 2, and to the second author and Ilias [EI2] in all
dimensions:
$$\lambda_1 ^c (M, [g]) \leq nV_c(M,[g])^{2/n}.$$

Our first result states that among all the possible conformal classes
of metrics on manifolds, the standard conformal class of the sphere
is the one having the lowest conformal spectrum.
\medskip

{\bf Theorem A}{ \it For any conformal class $[g]$ on $M$ and any
integer $k \geq 0$,}
$$\lambda_k ^c (M, [g]) \geq \lambda_k ^c ({\mathbb S}^n, [g_s]). $$

\medskip

Although the eigenvalues of a given Riemannian metric may have
nontrivial multiplicities, the conformal eigenvalues are all simple:
the conformal spectrum consists on a strictly
increasing sequence, and, moreover, the gap between two consecutive
conformal eigenvalues is uniformly bounded. Precisely, we have the
following theorem:

 \medskip

{\bf Theorem B}{ \it For any conformal class $[g]$ on $M$ and any
integer $k \geq 0$,
$$\lambda_{k+1} ^c (M, [g])^{n/2} - \lambda_k ^c (M, [g])^{n/2} \geq
\lambda_1 ^c
({\mathbb S}^n, [g_s]) = n^{n/2} \omega_n,$$
where $\omega_n $ is the volume of the $n$-dimensional Euclidean
sphere of radius one.}

\medskip

An immediate consequence of these two theorems is the following
explicit estimate of $ \lambda_k ^c (M, [g])$:

\medskip

{\bf Corollary 1} { \it For any conformal class $[g]$ on $M$ and any
integer }$k \geq 0$,
$$\lambda_k ^c (M, [g]) \geq n \omega_n^{2/n} k^{2/n}.$$

\medskip

Note that, for $k=1$, this last inequality has been recently proved
by Friedlander and Nadirashvili [FN] (see also [CE]). However, our
method is more general and simpler.

Of course, in the particular case of the $n$-sphere ${\mathbb S}^n$ endowed with its standard conformal class $[g_s]$, the equality holds in this inequality for $k=1$. The equality also holds for $k=2$ on ${\mathbb S}^2$ as it was recently proven by Nadirashvili [N2], i.e. $\lambda_2 ^c ({\mathbb S}^2,[g_s]) = 8\pi$.

Combined with the Korevaar estimate quoted above, Corollary 1 gives
$$ n \omega_n^{2/n} k^{2/n} \leq \lambda_k ^c (M, [g]) \leq C([g])
k^{2/n}.$$

Corollary 1 implies also that, if the $k$-th eigenvalue $\lambda_k
(g)$ of a metric $g$ is less than $n \omega_n^{2/n} k^{2/n}$, then
$g$ does not maximize $\lambda_k $ on its conformal class $[g]$. In particular, we have the following (negative) answer to a question of Yau concerning ${\mathbb S}^2$ (see [Y], p. 686):

\medskip

{\bf Corollary 2} {\it For any integer $k\ge 2$, the standard metric $g_s$ of ${\mathbb S}^2$ does not maximize $\lambda_k $, that is there exists a metric $g_k$ of volume one on ${\mathbb S}^2$ such that}
$$ \lambda_k (g_k) >\lambda_k (g_s).$$

\medskip
\noindent{Indeed, $\lambda_k (g_s) =4\pi [\sqrt{k}]([\sqrt{k}]+1)$, where $[\sqrt{k}]$ is the integer part of $\sqrt[]{k}$, while $\lambda_k ^c({\mathbb S}^2,[g_s]) \ge 8\pi k.$}
The same calculations show that, on ${\mathbb S}^3$, for any $ k\ge 2$, we have $\lambda_k (g_s)<3\omega_3^{2/3} k^{2/3}$, and then {\it the $k$-th eigenvalue does not achieve its maximum on $[g_s]$ at} $g_s$.

On the other hand, in any dimension we have, $\lambda_1^c({\mathbb S}^n,[g_s])=\lambda_1(g_s)=\lambda_2(g_s)=\cdots =\lambda_{n+1}(g_s)$. Consequently, $\forall k\in [2,n+1]$, {\it the standard metric $g_s$ of ${\mathbb S}^n$ does not maximize $\lambda_k$ in its conformal class}.

In [EI3] (see also [EI4] and [N]), Ilias and the second author
studied the property for a Riemannian metric to be critical (in a
generalized sense) for the functional $g \longmapsto \lambda_k(g) $.
A consequence of their results is that, if a metric $g$ is extremal
for $\lambda_k$ under conformal deformations, then the multiplicity
of $\lambda_k(g)$ is at least 2, which means that $\lambda_k(g)=
\lambda_{k+1}(g)$ or $\lambda_k(g)=\lambda_{k-1}(g)$. Combined with
Theorem B, this fact yields:

\medskip
{\bf Corollary 3} { \it If a Riemannian metric $g$ maximizes
$\lambda_1$ on its conformal class $[g]$, then it does not maximize
$\lambda_2$ on $[g]$.
More generally, a Riemannian metric $g$ cannot maximize
simultaneously three consecutive eigenvalues $\lambda_k$, $\lambda_
{k+1}$ and $\lambda_{k+2}$ on $[g]$. }

\medskip
 Applying Theorem B to an orientable surface $M_\gamma$ of genus
$\gamma$, we obtain the following result concerning the topological
spectrum.

\medskip

{\bf Corollary 4}{ \it For any fixed genus $\gamma$ and any integer
$k\geq 0$,
$$ \lambda_{k+1} ^{top} (\gamma) - \lambda_k ^{top} (\gamma) \geq 8
\pi,$$
and }
$$ \lambda_k ^{top} (\gamma) \geq 8(k-1)\pi+ \lambda_1 ^{top}
(\gamma) \geq 8k \pi.$$

\medskip

Our last result answers the following question: how does $ \lambda_k
^{top} (\gamma)$ behave as $\gamma$ increases?

\medskip

{\bf Theorem C} {\it For any fixed integer $k\geq0$, the function
$\gamma \longmapsto \lambda_k ^{top} (\gamma)$
is increasing, that is}

$$ \lambda_k ^{top} (\gamma +1) \geq \lambda_k ^{top} (\gamma).$$

\medskip

Recently, Brooks and Makover [BM] proved that, if $C$ is the Selberg
constant, then, for any $\varepsilon >0$, there exists an integer $N$
such that any compact orientable surface of genus $\gamma \geq N$
admits a hyperbolic metric $g$ with $ \lambda_1(g) \geq C -
\varepsilon$. As the area of such a hyperbolic surface is equal to $4
\pi (\gamma-1)$, it follows that
 $ \lambda_1 ^{top} (\gamma) \geq 4(C- \varepsilon) \pi (\gamma-1)$.
Although the Selberg conjecture "$C=1/4$" is still open, it has been
proved that $ C\geq 171/784 >1/5$ ([LRS]). Hence, for sufficiently
large $\gamma$, $ \lambda_1 ^{top} (\gamma) \geq {4\over 5} \pi
(\gamma-1)$ and then, $\forall k\geq 0$,
$$ \lambda_k ^{top} (\gamma) \geq {4\over 5} \pi (\gamma-1) + 8 \pi
(k-1).$$

\section{ Preliminary results}\label{2}

Roughly speaking, the proof of the two theorems A and B lies on the
following idea, also used in [FN]: locally, a Riemannian manifold $(M,g)$ is
almost Euclidean, and, consequently, almost conformal to the
sphere endowed with its standard metric $g_s$. Then, given a metric $h$
in the conformal class of the standard metric of the sphere, it will be possible to construct a
conformal deformation of $(M,g)$ around a point to make this
neighborhood arbitrarily close (in some sense) to $(\mathbb S^n,h)$.
To make these points precise, we will
establish some preliminary results and recall some facts from literature.

\medskip

Let $(M_{1},g_{1})$ and $(M_{2},g_{2})$ be two compact Riemannian
manifolds of the same dimension $n\geq 2$. Let us suppose that there
exists, for each $i\leq2$, a point $x_{i}\in M_{i}$ such that the
metric $g_i$ is flat in a neighborhood of $x_i$. Therefore, for
sufficiently small $\varepsilon >0$, the geodesic balls
$B_1(x_{1},\varepsilon) \subset M_{1}$ and
$B_2 (x_{2},\varepsilon) \subset M_{2}$ are both isometric to a
Euclidean ball. If $\Phi_{\varepsilon}: \partial B_1 (x_{1},
\varepsilon) \rightarrow
\partial B_2 (x_{2}, \varepsilon) $ is an induced isometry between
their boundaries, then we obtain a new closed manifold
$M_\varepsilon$ by glueing $M_1\backslash B_1(x_{1},\varepsilon)$ to
$M_2\backslash B_2(x_{2},\varepsilon)$ along $\Phi_{\varepsilon}$.
Let $\left\{ \lambda _k (\varepsilon);\hbox{\ } k\geq0 \right\}$ be
the spectrum of the natural Laplacian $\Delta_\varepsilon$ of
$M_\varepsilon$ (see [A]) associated to the piecewise smooth metric
$g_\varepsilon$ which coincide with $g_i$ on $M_i\backslash
B_i(x_{i},\varepsilon)$, and let $\left\{ \Lambda _k ;\hbox{\ }
k\geq0 \right\}$
be the reordered union of the spectra of $(M_{1},g_{1})$ and
$(M_{2},g_{2})$. Then we have the following:

\begin{lemma} For all $k\in {\mathbb N} $, we have
 $$\lim_{\varepsilon \to 0} \lambda_{k}(\varepsilon) =
\Lambda_{k}.$$
 \end{lemma}
This Lemma is a direct consequence of the min-max principle and the results
of [A]. Indeed, let $\{\mu_{k}(\varepsilon); \hbox{\ } k\geq0\}$ (resp.
$\{\nu_{k}(\varepsilon);\hbox{\ } k\geq0\}$) be the reordered union of
the spectra of $M_{1}\backslash B_1(x_{1}, \varepsilon) $
and $M_{2}\backslash B_2(x_{2}, \varepsilon) $ with the Dirichlet
(resp. Neumann) boundary condition. The following inequalities
are direct consequences of the min-max principle:
$$\nu_{k}(\varepsilon) \le \lambda_{k}(\varepsilon) \le
\mu_{k}(\varepsilon).$$
On the other hand, for each $i\leq2$, the spectrum of
$M_{i}\backslash B_i(x_{i}, \varepsilon) $
with the Dirichlet or the Neumann boundary condition,
converges, as $\varepsilon \to 0$, to the spectrum of the closed
manifold
$(M_{i},g_{i})$ (see [A]).

\bigskip
We also have a similar result in the case where we deform the metric on $M_{2}$
so that it collapses to a point. Indeed, by changing the scale of the metric $g_2$ if necessary, we may assume
that the radius one geodesic ball $B_2 (x_{2}, 1)$ of $(M_{2},g_{2})$
is contained in the flat neighborhood of $x_2$. Now, if we replace on
$M_2$ the metric $g_2$ by $g_{2}(\varepsilon) = \varepsilon^2 g_{2}$,
then the geodesic ball
 $B_1(x_{1}, \varepsilon) $ of $(M_{1},g_{1})$ becomes isometric to
the ball $B_2 (x_{2}, 1)$ of $(M_{2},g_{2})$ endowed with the new metric
$g_{2}(\varepsilon)$. Again, we consider the manifold $M_\varepsilon$
obtained by glueing $M_1\backslash B_1(x_{1},\varepsilon)$ to
$M_2\backslash B_2(x_{2},1)$ and endow it with the metric
$g_\varepsilon$ which coincide with $g_1$ on $M_{1}\backslash
B_1(x_{1}, \varepsilon) $ and with $g_{2}(\varepsilon)$ on
$M_{2}\backslash B_2 (x_{2}, 1)$. When $\varepsilon$ goes to zero,
Takahashi [T] proved that the spectrum of
$(M_{\varepsilon},g_{\varepsilon})$ converges to the spectrum of
$(M_{1},g_{1})$.

\begin{lemma} ([T]) For all $k\in {\mathbb N} $, we have
 $$\lim_{\varepsilon \to 0}
\lambda_{k}(M_{\varepsilon},g_{\varepsilon}) =
 \lambda_{k}(M_{1},g_{1}).$$
 \end{lemma}
Notice that, in dimension $n\geq 3$, this result can also be derived
from a theorem of Colin de Verdi\`ere [CV].

The proof of the theorems will also use the fact that the spectrum of
a Riemannian metric does not change much when we replace this latter
by a quasi-isometric one with a quasi-isometry ratio close to one.
Recall that two Riemannian metrics $g_1$ and $g_2$ on a compact
manifold $M$ are said to be $\alpha$-quasi-isometric, where
$\alpha\geq 1$, if, for any tangent vector $v\in TM$, $v \not = 0$,
 $$\frac{1}{\alpha ^2} \le \frac{g_{1}(v,v)}{g_{2}(v,v)} \le \alpha
^2.$$
The spectra of $g_1$ and $g_2$ are then related by the following
inequalities (see [D]): $\forall k \in {\mathbb N} ^*$,
$$\frac{1}{\alpha^{2(n+1)}} \le
\frac{\lambda_{k}(g_{1})}{\lambda_{k}(g_{2})} \le \alpha^{2(n+1)},$$
while their volumes satisfy
$$\frac{1}{\alpha ^n} \le \frac{V(g_{1})}{V(g_{2})} \le \alpha ^n.$$
The two following immediate observations will be useful in the sequel:

\begin {itemize}

\item[O1] If $g_1$ and $g_2$ are $\alpha$-quasi-isometric, then, for
any positive smooth function $f$ on $M$, the conformal metrics
$f^2g_1$ and $f^2g_2$ are also $\alpha$-quasi-isometric.

\item[O2] If $f_1$ and $f_2$ are two positive functions on $M$ such
that $\alpha^{-1} \le \frac{f_1(x)}{f_2(x)} \le \alpha$, then, for
any metric $g$ on $M$, the metrics $f_1^2g$ and $f_2^2g$ are
$\alpha$-quasi-isometric.

\end {itemize}

\begin{lemma} Let $(M,g)$ be a compact Riemannian manifold and let
$x_0$ be a point of $M$.

\smallskip
(i) For any positive $\delta$, there exists a
Riemannian metric $g_{_\delta}$ which is flat in a neighborhood of
$x_0$ and $(1+\delta )$-quasi-isometric to $g$ on $M$.

\smallskip
(ii) If, in addition, $g$ is conformally flat in a neighborhood of $x_0$,
then the metric $g_{_\delta}$ can also be chosen to be conformal to
$g$.
\end{lemma}

\noindent
\textit{Proof: } (i) In a normal coordinates system centered at $x_0$, we have
$$g_{ij}(x)=\delta _{ij} + O(\mid x \mid ^2).$$
Hence, it is clear that one can construct an adequate $g_{_\delta}$
by choosing it equal to $ \delta _{ij} $ in
a geodesic ball $B(x_0,r_\delta)$ of sufficiently small radius
$r_\delta$, and equal to $g$ in $M\backslash B(x_0, 2r_\delta)$.

\smallskip
(ii) In the case where $g$ is conformally flat in a neighborhood of
$x_0$,
we have the local expression
$$g_{ij}=f^2(x) \delta _{ij}, $$
where $f$ is a smooth function defined in a neighbourhood of $x_0$ and
such that $f(x_0)=1$. Thus, it suffices to take
$g_{_\delta}=\varphi_{_\delta} g$, where $\varphi_{_\delta}$ is
a positive smooth function on $M$ such that
$\varphi_{_\delta}=f^{-2}$ in a sufficiently small ball
$B(x_0,r_\delta)$, and $\varphi_{_\delta}=1 $ in $M\backslash B(x_0,
2r_\delta)$.

\section{ Proof of the theorems}\label{3}

\medskip
Let us start with the following elementary construction which will be useful in the sequel.

\paragraph {Construction:} Recall that the standard metric $g_s$ of the sphere ${\mathbb S}^n$ is expressed (via the
stereographic projection with respect to the north pole) by
$$g(x) = \frac{4}{(1+\Vert x \Vert^2)^2}g_{euc},$$
$g_{euc}$ being the Euclidean metric, while, given a positive number $R$, the metric
$$g_R(x)= \left \lbrace
\begin{array}{l}

 \frac{4}{(1+\Vert x \Vert^2)^2}g_{euc} \;\; if \; \Vert x \Vert \le R \\

\\

 \frac{4R^4}{(1+R^2)^2 \Vert x \Vert ^4}g_{euc} \;\; if \; \Vert x \Vert \ge R

\end{array}
\right.$$
corresponds to a metric on ${\mathbb S}^n$ which is flat in a ball around the north pole $N$, and coincide with $g_s$ outside this ball. The radius of this latter depends on $R$ and an easy
calculation shows that the set $\{\Vert x \Vert \ge R \}$, endowed with $g_R$, is
isometrically equivalent to an Euclidean ball $D_{\varepsilon(R)}$ of radius $\varepsilon(R) = \frac{2R}{(1+R^2)}$.

\smallskip
Now, for a metric $h$ conformal to $g_s$, we may consider a positive smooth
function $f$, with $f = 1$ on $\{\Vert x \Vert \ge R \}$, so that the metric $f^2(x) g_R$ represents $h$ outside the flat ball $D_{\varepsilon(R)}$. Moreover, up to a scaling of the variable $x$, it is possible to prescribe the radius of the flat ball. Indeed, it suffices to consider, for any positive $\rho$, the metric
$$g_{f,R, \rho}(x) = f^2(\frac{Rx}{\rho}) g_R({Rx\over \rho})=\left \lbrace
\begin{array}{l}
f^2(\frac{Rx}{\rho})\frac{R^2}{\rho^2}
\frac{4}{(1+ \frac{R^2}{\rho^2}\Vert x \Vert^2)^2}g_{euc} \;\; if \; \Vert x \Vert
\le \rho \\

\\

\frac{4R^2\rho^2}{(1+R^2)^2 \Vert x \Vert ^4}g_{euc} \; \; if
\; \Vert x \Vert \ge \rho

\end{array}
\right.$$
Note that if $\Vert x \Vert = \rho$, then the metric at $x$
is $\frac{\varepsilon^2(R)}{\rho^2}g_{euc}$.

\medskip

Before going further into the proof, let us note that the
conformal eigenvalues are not necessarily achieved by smooth metrics, so that we will always work with
smooth Riemannian metrics whose eigenvalues are almost extremal for the
problem we study, and then pass to the limit.

\paragraph {Proof of Theorem A:} Let $(M,g)$ be a $n$-dimensional compact Riemannian manifold of
volume 1 and let $k$ be a positive integer.
Let us fix a positive real number $\delta$.

\smallskip
On the sphere $\mathbb S^n$, consider a Riemannian
metric $h\in [g_s]$ of volume one satisfying
$$\lambda_k(h) \geq \lambda_{k}^c(\mathbb S^n,[g_{s}])- \delta.$$
As it belongs to the standard conformal class $[g_s]$, the metric $h$
is locally conformally flat and, applying Lemma 2.3 (ii), there exists a
metric $h_\delta$, conformal to $g_s$, flat in the geodesic ball
$D_{r_\delta} \subset {\mathbb S}^n$ of radius $r_{\delta}$
and $(1+\delta)$-quasi-isometric
to $h$. In particular, $h_\delta$ satisfies

\begin{eqnarray}
\nonumber {} \lambda_k(h_\delta)V(h_\delta)^{2/n} &\ge & (1+\delta)^{-2(n+2)}\lambda_k(h)V(h)^{2/n} \\
\nonumber {} &\ge & (1+\delta)^{-2(n+2)}(\lambda_{k}^c(\mathbb S^n,[g_{s}])- \delta).
\end{eqnarray}

On the other hand, let $g_{\delta}$ be a metric on $M$ satisfying
the conditions of Lemma 2.3 (i). We multiply $g_{_\delta}$ by a constant
$C_\delta^2$ so that the geodesic ball $B (x_0,2)$ becomes flat. Now,
as explained in the construction above,
for any positive $\varepsilon < r_\delta$, there exists on $M$ a
 metric $f_{\varepsilon}^2g_{\delta}$ conformal to $C_{\delta}^2g_{\delta}$
 such that

 \begin{itemize}
 \item[-] The closed ball $\bar B(x_0,1)\subset M$ endowed with $f_{\varepsilon}^2g_{\delta}$
becomes isometric to $({\mathbb S}^n \backslash D_\varepsilon, h_\delta)$.

\item[-] The metric $f_{\varepsilon}^2g_{\delta}$ coincide with $\varepsilon^2 C_{\delta}^2g_{\delta}$ on $M \backslash B(x_0,1)$.
\end{itemize}
Hence, we may
identify the manifold $(M,f_{\varepsilon}^2g_{\delta})$
to the manifold $M_\varepsilon$ of Lemma
2.2 above obtained by glueing $M \backslash B(x_0,1)$ to ${\mathbb
S}^n \backslash D_\varepsilon$. Lemma 2.2 tells us that
$\lambda_k(g_\varepsilon)$ converges, as $\varepsilon \to 0$, to
$\lambda_k(h_\delta)$. It is also clear that the volume
$V(g_\varepsilon)$ of $(M,g_\varepsilon)$ converges to the volume of
$({\mathbb S}^n , h_\delta)$. Therefore, there exists $\varepsilon
>0$ such that
$$ \lambda_k(g_\varepsilon) V(g_\varepsilon)^{2/n}\geq
\lambda_k(h_\delta)V(h_\delta)^{2/n} -\delta \geq
(1+\delta)^{-2(n+2)}(\lambda_{k}^c({\mathbb S}^n,[g_{s}])-
\delta)-\delta.$$

Now, using classical density results and the observation O2 above, we
may find a smooth function $\bar f_\varepsilon$ on $M$ so that the
smooth metric $\bar g_\varepsilon=\bar f_\varepsilon ^2 g_{_\delta}$
is $(1+\delta)$-quasi-isometric to $g_\varepsilon=f_\varepsilon^2
g_{_\delta}$. As $g_{_\delta}$ is $(1+\delta)$-quasi-isometric to
$g$, the observation O1 above tells us that the metric
$g'_\varepsilon = \bar f_\varepsilon ^2 g$ is in fact
$(1+\delta)^2$-quasi-isometric to $g_\varepsilon$. Therefore, we have
\begin{eqnarray}
\nonumber {} \lambda_k(g'_\varepsilon) V(g'_\varepsilon)^{2/n} &\geq&
(1+\delta)^{-4(n+2)}\lambda_k(g_\varepsilon) V(g_\varepsilon)^{2/n}\\
 \nonumber{} &\ge& (1+\delta)^{-6(n+2)}(\lambda_{k}^c({\mathbb
S}^n,[g_{s}]) -\delta) -\delta (1+\delta)^{-4(n+2)} \\
\nonumber{} &=& \lambda_{k}^c({\mathbb S}^n,[g_{s}]) - O(\delta),
\end{eqnarray}
with $O(\delta) \to 0$ as $\delta \to 0$. Since $g'_\varepsilon$ is
conformal to $g$, it follows, according the definition of
$\lambda_{k}^c (M, [g])$,
$$\lambda_{k}^c (M, [g]) \ge \lambda_{k}^c({\mathbb S}^n,[g_{s}]) -
O(\delta).$$
As $\delta$ can be chosen arbitrarily small, we get the desired
inequality:
$$\lambda_{k}^c (M, [g]) \ge\lambda_{k}^c({\mathbb S}^n,[g_{s}]).$$
\hbox{} \hfill $\Box$

\paragraph {Proof of Theorem B:} Let $(M,g)$ be a $n$-dimensional compact Riemannian manifold of
volume 1. Let $k$ be a positive integer and $\rho$ a positive real
number. We endow ${\mathbb S}^n$ with the metric
$$h_\rho={n {\omega_n ^{2/n}} \over \lambda_{k}^c (M, [g])+\rho } \;
g_s,$$
so that its first positive eigenvalue becomes equal to $
\lambda_{k}^c (M, [g])+\rho$ (recall that $\lambda_1(g_s)=n {\omega_n
^{2/n}}$), and consider a metric $g_\rho \in [g]$ on $M$ such that
$$ \lambda_k (g_\rho) \ge \lambda_{k}^c (M, [g]) - \rho/2 .$$

Let $\delta$ be a sufficiently small positive real number so that

$$(1+\delta)^{2(n+1)} \lambda_{k}^c (M, [g]) \le (1+\delta)^{-2(n+1)
} (\lambda_{k}^c (M, [g])+\rho).$$
 We apply Lemma 2.3 to get a metric $g_{\rho,\delta}$ on $ M$, a
metric $h_{\rho,\delta}\in [g_s]$ on $ {\mathbb S}^n$, and a constant
$r_\delta >0$ such that

\begin {itemize}

\item[-] $g_{\rho,\delta}$ is $(1+\delta)$-quasi-isometric to
$g_\rho$ and $h_{\rho,\delta}$ is $(1+\delta)$-quasi-isometric to
$h_\rho$,

\item[-] $\forall \; \varepsilon \in (0,r_\delta)$, the geodesic
balls $B(x_0,\varepsilon) \subset (M,g_{\rho,\delta})$ and
$D(x_1,\varepsilon) \subset ({\mathbb S}^n,h_{\rho,\delta})$ are
isometric to a Euclidean ball, where $x_0$ and $x_1$ are two given
points of $M$ and ${\mathbb S}^n $ respectively.

\end {itemize}

\noindent
As in the proof of Theorem A, we notice that, for any $\varepsilon
\in (0,r_\delta)$, the closed ball $(\bar
B(x_0,\varepsilon),g_{\rho,\delta}) $ is conformally equivalent to
$({\mathbb S}^n \backslash D(x_1,\varepsilon), h_{\rho,\delta})$, and
that there exists a piecewise smooth function $f_\varepsilon$ on $M$
which is equal to $1$ on $M \backslash B(x_0,\varepsilon)$ and such
that $(\bar B (x_0,\varepsilon), f_\varepsilon ^2 g_{\rho,\delta})$
is isometric to $({\mathbb S}^n \backslash D(x_1,\varepsilon),
h_{\rho,\delta})$. Now, the manifold $(M,g_\varepsilon=f_\varepsilon
^2g_{\rho,\delta})$ is identified to the manifold $M_\varepsilon$ of
Lemma 2.1 obtained by glueing $(M \backslash B(x_0,\varepsilon),
g_{\rho,\delta})$ to $({\mathbb S}^n \backslash D(x_1,\varepsilon),
h_{\rho,\delta}) $. According to this lemma, the spectrum of
$(M,g_\varepsilon)$ converges, as $\varepsilon$ goes to $0$, to the
reordered union of the spectra of $(M,g_{\rho,\delta})$ and
$({\mathbb S}^n, h_{\rho,\delta})$. From the construction, we have
the following inequalities:
$$ \lambda_k (g_{\rho,\delta})\le (1+\delta)^{2(n+1)}
\lambda_{k}(g_\rho) \le (1+\delta)^{2(n+1)} \lambda_{k}^c (M, [g]),$$
and
$$ \lambda_1 (h_{\rho,\delta})\ge (1+\delta)^{-2(n+1) } (\lambda_1
(h_\rho)\ge(1+\delta)^{-2(n+1) } (\lambda_{k}^c (M, [g])+\rho).$$
Hence, from the smallness condition above satisfied by $ \delta$, we
have
$$\lambda_k (g_{\rho,\delta})\le\lambda_1 (h_{\rho,\delta}).$$
Consequently, the lowest $(k+2)$ eigenvalues in
$Spec(g_{\rho,\delta}) \cup Spec (h_{\rho,\delta})$ are:
$$\lambda_0(h_{\rho,\delta}),\; \lambda_0(g_{\rho,\delta}),\;
\lambda_1(g_{\rho,\delta}), \lambda_2(g_{\rho,\delta}), \cdots
\lambda_k(g_{\rho,\delta}).$$
Thus,
$$\lim_{\varepsilon \to 0}
\lambda_{k+1}(g_\varepsilon)=\lambda_k(g_{\rho,\delta}).$$

On the other hand, the volume of $(M,g_\varepsilon)$ converges, as
$\varepsilon \to 0$, to
$V(g_{\rho,\delta}) + V(h_{\rho,\delta})$.
Since $g_{\rho,\delta}$ and $h_{\rho,\delta}$ are
$(1+\delta)$-quasi-isometric to $g_\rho$ and $h_\rho$ respectively,
we have
$$\lambda_k(g_{\rho,\delta})\ge (1+\delta)^{-2(n+1)}
\lambda_{k}(g_\rho)\ge (1+\delta)^{-2(n+1) } (\lambda_{k}^c (M,
[g])-\rho /2),$$
and

\begin{eqnarray}
\nonumber {} V(g_{\rho,\delta}) + V(h_{\rho,\delta}) &\ge& (1+\delta)^{-n} (V(g_\rho)+V(h_\rho)) \\
 \nonumber{} &=& (1+\delta)^{-n} (1+{n^{n/2} {\omega_n } \over
(\lambda_{k}^c (M, [g])+\rho)^{n/2} }).
\end{eqnarray}

Therefore, there exists a positive $\varepsilon$ such that
\begin{eqnarray}
\nonumber {} \lambda_{k+1}(g_\varepsilon)^{n/2} V(g_\varepsilon) &\ge& (1+\delta)^{-n(n+2) }\times\\
\nonumber{} &{}&\times (\lambda_{k}^c (M, [g])-\rho )^{n/2} (1+{n^{n/2} {\omega_n } \over (\lambda_{k}^c (M, [g])+\rho)^{n/2}}).
 \end {eqnarray}
As in the proof of Theorem A, we use the observations O1 and O2 to
get a smooth metric $\bar g_\varepsilon$ conformal to $g$ and
$(1+\delta)^2$-quasi-isometric to $g_\varepsilon$. Hence,
$$ \lambda_{k+1}^c (M, [g])^{n/2} \ge \lambda_{k+1}(\bar
g_\varepsilon)^{n/2} V(\bar g_\varepsilon) \geq
(1+\delta)^{-2n(n+2)}\lambda_{k+1}( g_\varepsilon)^{n/2} V(
g_\varepsilon),$$
and then,
$$ \lambda_{k+1}^c (M, [g])^{n/2} \ge$$
$$\ge (1+\delta)^{-3n(n+2) }
(\lambda_{k}^c (M, [g])-\rho )^{n/2} (1+{n^{n/2} {\omega_n } \over
(\lambda_{k}^c (M, [g])+\rho)^{n/2} }),$$
which gives, as $\delta \to 0$ and, then, $\rho \to 0$,
$$ \lambda_{k+1}^c (M, [g])^{n/2}\ge \lambda_{k}^c (M, [g])^{n/2} +
n^{n/2} \omega_n .$$ \hbox{} \hfill $\Box$

\paragraph {Proof of Theorem C:} Let $M\gamma$ be a compact
orientable surface of genus $\gamma$ and let $k$ be a natural
integer. Given a positive real number $\delta$, let $g_{_\delta}$ be
a Riemannian metric of area one on $M\gamma$ such that
$$\lambda_k (g_{_\delta}) > \lambda_k ^{top} (\gamma) - \delta/2.$$

Let us attach to $M_\gamma$ a "thin" handle of radius $\varepsilon >
0$ and
length $l > 0$ as described in [A]. After smoothing, we get a
compact surface $M_{\gamma +1}$ of genus $(\gamma + 1)$ endowed with
a Riemannian
metric $ g_{\delta , \varepsilon , l}$ such that, when $\varepsilon$
goes to
zero,
\begin {itemize}

 \item[-] the spectrum of the Laplacian of $g_{\delta , \varepsilon ,
l}$ converges to the reordered union of the spectrum of
$(M_\gamma, g_{_\delta})$ and the spectrum of the segment $[0,l]$ for
the Laplacian
with Dirichlet boundary condition.
\item[-] the area $V (g_{\delta , \varepsilon , l}) $ of $g_{\delta ,
\varepsilon , l}$ converges to $1$,

\end {itemize}

Choosing $l$ sufficiently small, one can suppose that the first
Dirichlet eigenvalue of $[0,l]$ is greater than $\lambda_k ^{top}
(\gamma)$. Hence, the lowest $(k+1)$ eigenvalues in $Spec(M_\gamma,
g_{_\delta}) \cup Spec ([0,l])$ are:
$$ \lambda_0(g_{_\delta}),\; \lambda_1(g_{_\delta}),
\lambda_2(g_{_\delta}), \cdots \lambda_k(g_{_\delta}),$$
and then
$$\lim_{\varepsilon \to 0}\lambda_k (g_{\delta , \varepsilon , l}) =
\lambda_k(g_{_\delta}).$$
 Therefore, there exist two positive constants $\varepsilon$
and $l$ such that
$$\lambda_k (g_{\delta , \varepsilon , l}) \geq \lambda_k
(g_{_\delta}) - \delta/2 >\lambda_k ^{top} (\gamma) -
\delta$$
 and
$$V (g_{\delta , \varepsilon , l}) > 1 - \delta.$$
Consequently

$$\lambda_k ^{top}(\gamma + 1) \geq \lambda_k (g_{\delta ,
\varepsilon , l}) V (g_{\delta , \varepsilon , l})
\geq (\lambda_k ^{top} (\gamma) - \delta) (1-\delta).$$
In conclusion, we have
$$\lambda_k ^{top} (\gamma + 1) \geq \lambda_k ^{top} (\gamma).$$
\hbox{} \hfill $\Box$

\end{document}